\documentstyle{amsart}

\numberwithin{equation}{section}

\newtheorem{theorem}{Theorem}[section]

\newtheorem{corollary}[theorem]{Corollary}
\newtheorem{proposition}[theorem]{Proposition}

\theoremstyle{remark}
\newtheorem{remark}[theorem]{Remark}

\newtheorem{definition}[theorem]{Definition}
\newtheorem{example}[theorem]{Example}          
\newtheorem{ack}{Acknowledgment}

\newcommand{\Lie}{\operatorname{Lie}}

\begin{document}

\title{Twistor quotients of hyperk\"ahler manifolds}
\author{Roger Bielawski}
\thanks{Research supported by an EPSRC Advanced Research Fellowship}

\address{Department of Mathematics\\
University of Glasgow\\
Glasgow G12 8QW\\
Scotland  }

\email{rb@@maths.gla.ac.uk}

\begin{abstract} We generalize the hyperk\"ahler quotient construction to the situation where there is no group action preserving the hyperk\"ahler structure but for each complex structure there is an action of a complex group preserving the corresponding complex symplectic structure. Many (known and new) hyperk\"ahler manifolds arise as quotients in this setting. For example, all hyperk\"ahler structures on semisimple coadjoint orbits of a complex semisimple Lie group $G$ arise as such quotients of $T^\ast G$. The generalized Legendre transform construction of Lindstr\"om and Ro\v{c}ek is also explained in this framework.     \end{abstract}

\maketitle

\begin{center} \begin{sc} Introduction\end{sc}\end{center}

The motivation for this work stems from two problems. The first is the following question: when is a complex-symplectic quotient of a hyperk\"ahler manifold hyperk\"ahler?  A good example is the hyperk\"ahler structure on  $M=T^\ast G$, where $G$ is a complex semisimple Lie group (found by Kronheimer, cf. \cite{BielJLMS}). The complex symplectic quotients of $M$ by $G$ are precisely coadjoint orbits of $G$. These carry hyperk\"ahler structures by the work of Kronheimer \cite{Kr}, Biquard \cite{Bi} and Kovalev \cite{Kov}.
\par
The second motivating problem is the generalized Legendre transform (GLT) construction of hyperk\"ahler metrics due to Lindstr\"om and Ro\v{c}ek \cite{LR}. Unlike the ordinary Legendre transform which produces $4n$-dimensional hyperk\"ahler metrics with $n$ commuting Killing vector fields, the GLT produces metrics without (usually) any Killing vector fields. The defining feature of these metrics is that their twistor space admits a holomorphic projection onto a vector bundle of rank $n$ over ${\Bbb C} P^1$. 
\par
It turns out that in both of these problems there is a group-like object involved, namely a bundle of complex groups over ${\Bbb C} P^1$ which act fiberwise on the twistor space $Z$ of a hyperk\"ahler manifold. This action is also Hamiltonian for the twisted symplectic form of $Z$. Thus, whenever we have such an action, we can form fiberwise complex-symplectic quotient of $Z$ giving us (in good cases) a new twistor space. Similarly, in the case of the GLT, the projection onto the vector bundle $V$ should be regarded as the moment map for an action of a bundle of abelian groups on $Z$, which preserve the twisted symplectic form.
\par
We call our bundles of groups over ${\Bbb C}P^1$ {\it twistor groups}. The simplest definition of a twistor group is a {\it group in the category of spaces over ${\Bbb C}P^1$ with a real structure}.
\par
If a twistor group acts on the twistor space $Z$ of a hyperk\"ahler manifold $M$, we can interpret (in most cases) the resulting vector fields on $Z$ as objects on $M$, namely either as higher rank Killing spinors (cf. \cite{Dun}) or, in the $E-H$ formalism (cf. \cite{Sal}) as sections of $E\otimes S^{2i+1}H$ ($i>1$) satisfying equations analogous to the Killing vector field equation (case $i=1$).
\par
The main purpose of this paper is to introduce the concept of twistor groups and their actions and to give some interesting examples. We also prove results which can be viewed as new constructions of hyperk\"ahler manifolds.\newline

\section{Twistor groups and their actions\label{Action}}

\subsection{Twistor groups} 

Let $X$ be a complex manifold. A {\em space over $X$} is a complex space $Z$ together with a surjective holomorphic map ({\em projection}) $\pi:Z\rightarrow X$. We shall say that $Z@>\pi>> X$ is {\em smooth} if $Z$ is smooth and $\pi$ is a submersion.
\par
The category of spaces over $X$ is a category with products (fiber product) and a final object ($X@>{\rm Id} >> X$). In any category with such properties we can define a {\em group} as an object ${\cal G}$ together with morphisms defining group multiplication, inverse, and the identity. Thus we define:

\begin{definition} A {\em group over $X$} is a group in the category of spaces over $X$.\label{group}\end{definition}

More explicitly, a group over $X$ is a space ${\cal G}@>\pi>> X$ together with fibrewise   holomorphic maps  $\;\cdot:{\cal G}\times_X {\cal G}\rightarrow {\cal G}$ ({\em multiplication}), $i:{\cal G}\mapsto {\cal G}$ ({\em group inverse}) and $1:X\rightarrow {\cal G}$ ({\em identity section}) which commute with $\pi$ and satisfy the group axioms. In particular, for each  $x\in X$ $\bigl(\pi^{-1}(x),\,\cdot\,,i_{|_{\pi^{-1}(x)}},1(x)\bigr)$ is a group.

\begin{remark} Even if one is interested (as we are) primarily in smooth groups over $X$, one cannot avoid the singular ones, since a subgroup of a smooth group  can be singular. In particular the stabilizers of smooth group actions can be singular.\end{remark}

We shall be interested mostly in the case when $X={\Bbb C}P^1$ and the spaces over ${\Bbb C}P^1$ come equipped with an antiholomorphic involution ({\em real structure}) covering the antipodal map on ${\Bbb C}P^1$. The category of spaces with a real structure over ${\Bbb C}P^1$  is also a category with products and a final object. Therefore we can define:

\begin{definition} A {\em twistor group} is a  group in the category of smooth spaces  with a real structure over ${\Bbb C}P^1$.\end{definition}

\begin{remark} Although the natural setting is the category of complex spaces rather than of manifolds, all our examples involve only smooth groups. In addition, the proofs are either simpler or work only for smooth groups. \end{remark}

Let us give a few examples of  twistor groups.
\begin{example} Let $G$ be a complex Lie group equipped with an antiholomorphic involutive automorphism $\sigma$. Then $G\times {\Bbb P}^1$ with the involution $(\sigma, a)$, where $a$ is the antipodal map, is a twistor group which we shall call a {\em trivial twistor group} (with structure group $G$) and denote by $\underline{G}$.\label{trivial}\end{example}
\begin{example} Many nontrivial examples arise as twistor subgroups of $\underline{G}$. For example, if $G$ acts fibrewise on a  space $Z$ with a real structure over ${\Bbb C}P^1$, then the stabilizer of any real section of $Z$ is a twistor subgroup of $\underline{G}$. In particular, we can take the adjoint action of $G$ on $Z=\underline{\frak g}\otimes V$, where $V$ is a vector bundle over ${\Bbb C}P^1$ equipped with a real structure.label{sub}\end{example}  

\begin{example} Another important twistor group is constructed as follows.  Let $G$ be reductive and let ${\frak  k}$ denote the Lie algebra of the maximal compact subgroup of $G$. Let $\rho:{\frak su}(2)\rightarrow {\frak  k}$ be a homomorphism of Lie algebras.
For each element $z=(a,b,c)$, $a^2+b^2+c^2=1$, of $S^2\simeq {\Bbb C}P^1$, define a subalgebra ${\frak  n}_z$ of ${\frak  g}$ as the sum of negative eigenspaces of $\text{ad}\bigl(a\rho(\sigma_1)+b\rho(\sigma_2)+c\rho(\sigma_3)\bigr)$, where $\sigma_1,\sigma_2,\sigma_3$ denote the Pauli matrices. Now define ${\cal N}$ as a twistor subgroup of $\underline{G}$ whose fiber at $z$ is the subgroup of $G$ whose Lie algebra is $N_z$. It is straightforward to observe that the real structure of $\underline G$ acts on ${\cal N}$. We also observe that each fiber of ${\cal N}$ is the unipotent radical of the parabolic subgroup of $G$ determined by $\rho$. \label{borel}\end{example}
\begin{example} A vector bundle over ${\Bbb P}^1$ equipped with a (linear) real structure is an abelian twistor group. We observe that a line bundle ${\cal O}(k)$ is a twistor group if and only if $k$ is even.\end{example}

The last example can be generalized as follows. Let ${\cal G}$ be any twistor group. For an open subset $U$ of ${\Bbb C}P^1$ denote by ${\cal G}_U$ the group over $U$ obtained as the restriction of ${\cal G}$ to $U$. Now suppose that we are given a covering $\{U_i\}$ of ${\Bbb C}P^1$, invariant under the antipodal map and a fibrewise automorphism $\{\phi_{ij}\}$ of ${\cal G}_{U_i\cap U_j}$ for each nonempty intersection $U_i\cap U_j$. In addition we suppose that the family of $\phi_{ij}$  is $\tau$-equivariant, where $\tau$ is the real structure of ${\cal G}$. Then gluing together ${\cal G}_{U_i}$ via the $\phi_{ij}$ gives us a new twistor group, locally isomorphic to ${\cal G}$. We deduce the following:
\begin{proposition} Let ${\cal G}$ be a twistor group. Then the isomorphism classes of twistor groups locally isomorphic to ${\cal G}$ are in bijective correspondence with elements of (non-abelian) sheaf cohomology group $H^1_{\Bbb R}\bigl({\Bbb C}P^1,{\cal A}\bigr)$, where ${\cal A}(U)$ is the group of automorphisms of ${\cal G}_U$.\hfill$\Box$\label{aut}\end{proposition}
The subscript ${\Bbb R}$ denotes $\tau$-invariant elements.
\par
In particular, if $G$ is a complex Lie group with an antiholomorphic automorphism, then we can consider twistor groups which are locally isomorphic to $\underline{G}$. We shall call such twistor groups {\em locally trivial}. We have:
\begin{corollary} The isomorphism classes of locally trivial twistor groups with structure group $G$ are in 1-1 correspondence with elements of $H^1_{\Bbb R}\bigl({\Bbb C}P^1,{\cal O}(\text{Aut}(G))\bigr)$. \hfill$\Box$\end{corollary}

We shall call a twistor group ${\cal G}$ {\em discrete}, if each fiber of ${\cal G}$ is discrete. The following example shows that twistor groups need not be locally trivial.
\begin{example} Let $s$ be a real meromorphic section of ${\cal O}(-2)$ having poles at a pair of antipodal points $a,-1/\bar{a}$. We have a (smooth) discrete twistor subgroup ${\cal L}$ of ${\cal O}(-2)$ defined as the subgroup of ${\cal O}(-2)$ with fiber $0$ at $a$ and $-1/\bar{a}$ and ${\Bbb Z}s(x)$ at other points. \label{L}\end{example}

Notice that ${\cal O}(-2)/{\cal L}$ (fibrewise quotient) is also a twistor group.

\subsection{Twistor Lie algebras} Let ${\cal G}@>\pi>> X$ be a smooth group over $X$. Then the normal bundle to the identity section has a natural structure of a {\em Lie algebra over $X$} (i.e. a Lie algebra in the category of vector bundles over $X$). We shall denote this space by $\Lie({\cal G})$. In particular $\Lie({\cal G})$ is locally trivial as a vector bundle (cf. \cite{EMS}). 

Let us consider the structure of a twistor Lie algebra ${\cal L}$ in a more detail. As observed in the previous section, ${\cal L}$ is a locally trivial vector bundle and so it splits as $\bigoplus {\cal O}(p_i)$ for some integers $p_1,\dots, p_n$. We choose coordinates $e_1,\dots,e_n$ for ${\cal L}$ over $\zeta\neq \infty$ and $\tilde{e}_1,\dots,\tilde{e}_n$ over $\zeta\neq 0$, so that $\tilde{e}_i=\zeta^{-p_i}e_i$ over ${\Bbb C}P^1-\{0,\infty\}$. The fibrewise Lie bracket is given by
\begin{equation} \bigl[e_i,e_j\bigr]_\zeta=\sum_{k=1}^n f_k(\zeta)e_k,\label{fk}\end{equation}
\begin{equation} \bigl[\tilde{e}_i,\tilde{e}_j\bigr]_{\tilde{\zeta}}=\sum_{k=1}^n \tilde{f}_k(\tilde{\zeta})\tilde{e}_k,\label{tildefk}\end{equation}
for some holomorphic functions $f_k=f_k^{ij}$, $\tilde{f}_k=\tilde{f}^{ij}_k$.
Comparing the two expressions for the bracket over $\zeta\neq 0,\infty$, we see that $f_k,\tilde{f}_k$ define a section of ${\cal O}(p_i+p_j-p_k)$. In particular if some of these sections are nonconstant (i.e. have zeros), then ${\cal L}$ is locally nontrivial as a bundle of Lie algebras.

The preceding considerations imply the following fact, which will be useful later on.
\begin{proposition} Let ${\cal G}$ be a smooth twistor group whose Lie algebra splits as a sum of line bundles of negative degrees. Then ${\cal G}$ is nilpotent.\hfill $\Box$ \label{neg}\end{proposition}

\subsection{Actions of twistor groups} 
We now define actions of twistor groups or of groups over $X$. Once more, this is a tautological definition in any category with products and a final object. In our case  an {\em action} of a group ${\cal G}@>\pi>> X$ on $Z@>p>> X$ is a holomorphic map $$ \cdot:{\cal G}\times_X Z\rightarrow Z$$
which commutes with the projections and which is a group action on each fiber. An action of a twistor group is required to respect the real structures ${\cal G}$ and $Z$ (i.e. $\tau(g\cdot z)=\tau(g)\cdot \tau(z)$).
Most notions related to (ordinary) group actions carry over to actions of groups over $X$. Thus, we shall say that the action of ${\cal G}$ is free (resp. locally free) if each fiber action is free (resp. locally free). We can define equivariant morphisms. We also observe that for smooth groups over $X$ we have canonical notions of the {\em adjoint} and {\em coadjoint} action. Finally we can define orbits and stabilizers of sections of $Z@>p >>X$.

\begin{remark} In the case of an action on a twistor space of a hyperk\"ahler manifold $M$, we shall also speak of ${\cal G}$ acting on $M$. Similarily, if $s$ is a twistor line corresponding to a point $m$ in $M$, we can speak of the stabilizer of ${\cal G}$ at $m\in M$ etc.\end{remark}

We shall be particularly interested in the following types of actions.

\begin{definition} Let a smooth twistor group ${\cal G}@>\pi>> {\Bbb C}P^1$ act on a twistor space $Z@>p>> {\Bbb C}P^1$. We shall say that the action is {\em symplectic} (resp. {\em Hamiltonian}), if the action is symplectic on each fiber for the twisted symplectic form $\omega$ on $Z$ (resp. if it is symplectic and if there is a holomorphic map $\mu:Z\rightarrow \text{Lie}({\cal G})^\ast\otimes {\cal O}(2)$ which is the moment map for the twisted symplectic form $\omega$ on each fiber). \end{definition}

\begin{example} Let a compact Lie group $K$ act on a hyperk\"ahler manifold $M$ by hyperk\"ahler isometries and suppose that this action extends to the action of $K^{\Bbb C}$ for each complex structure of $M$. Then the trivial group $\underline{K^{\Bbb C}}$ acts symplectically on $Z$. Consequently, any twistor subgroup of $\underline{K^{\Bbb C}}$ acts symplectically on $Z$. If the action of $K$ is tri-Hamiltonian, then the action of $\underline{K^{\Bbb C}}$ and any of its twistor group subspaces is Hamiltonian.\end{example} 
\begin{example} Let $M={\Bbb H}$ so that $Z={\cal O}(1)\oplus {\cal O}(1)$. Then the twistor group ${\cal O}(1)\oplus {\cal O}(1)$ acts on $Z$ via fibrewise addition. This action is symplectic, but not Hamiltonian.\end{example}
\begin{example} ({\bf Atiyah-Hitchin manifold}). This is an example of a Hamiltonian action of a twistor group (namely ${\cal O}(-2)$) which does not arise from any hyperk\"ahler group action. Let $M$ be the hyperk\"ahler manifold of strongly centered monopoles of charge $2$, i.e. the double (or universal) cover of the Atiyah-Hitchin manifold. With respect to any complex structure it is biholomorphic to the space of degree $2$ rational maps of the form $p(z)/q(z)$ where $q(z)=z^2-c$ and $p(z)=az+b$ with $b^2-c a^2=1$. Let $\beta$ denote any of the two roots of $q$. The complex symplectic form is given (on the set where $c\neq 0$) by $\omega=\frac{dp(\beta)}{p(\beta)}\wedge d\beta$. The twistor space $Z$ of $M$ is essentially given by requiring that $\beta$ is a section of ${\cal O}(2)$ while $p(\beta)$ is a value of a certain line bundle $L^{-2}$ on ${\cal O}(2)$ \cite{AH}. We define an action of ${\Bbb C}$ on $M$ by
$$ \lambda\cdot\begin{pmatrix} a\\b\\c\end{pmatrix}=\begin{pmatrix} \cosh\lambda\beta & \frac{\sinh\lambda\beta}{\beta} & 0\\
\beta\sinh\lambda\beta & \cosh\lambda\beta & 0\\ 0 & 0 & 1\end{pmatrix} \begin{pmatrix} a\\b\\c\end{pmatrix}.  $$
Here $\beta=\pm \sqrt{c}$. This action sends $p(\beta)$ to $e^{\lambda\beta}p(\beta)$ and so it respects the complex symplectic form $\omega$. By looking at the transition functions for the twistor space $Z$ we conclude that this action extends to the fibrewise action of ${\cal O}(-2)$ on $Z$. One can check that the real structures are compatible with this action.  This action is locally free and the orbits of twistor lines are isomorphic to the twistor groups $O(-2)/{\cal L}$ defined in example \ref{L}. Finally, this action is Hamiltonian and the twistor  lines of $Z$ project via the moment map $\mu:Z\rightarrow {\cal O}(4)$ to the {\em spectral curves} of the monopoles. \end{example}

We have an obvious restriction on twistor groups which can act locally freely at any twistor line.
\begin{proposition} Suppose that a twistor group ${\cal G}$ acts on a twistor space $Z$ and that the action is locally free at some twistor line $s$. Then $\Lie{\cal G}$ is the sum of line bundles of degree at most one.  \label{restriction}\end{proposition}
\begin{pf} We have an injective morphism $\imath:{\cal L}\rightarrow {\cal O}(1)\otimes {\Bbb C}^n$ of vector bundles over ${\Bbb C}P^1$. \end{pf}
For Hamiltonian actions the restriction is more severe.
\begin{proposition} Suppose that there is a Hamiltonian action of a twistor group ${\cal G}$ on a twistor space $Z$ which is  locally free at some twistor line. Then $\Lie{\cal G}$   is the sum of line bundles of degree at most zero.\hfill $\Box$\label{restrict}\end{proposition}

\subsection{Quotients and principal bundles}

An action of a twistor group ${\cal G}$ on a space $Z@>p>> {\Bbb C}P^1$ defines an equivalence relation $R\subset Z\times Z$: $(z_1,z_2)\in R\iff p(z_1)=p(z_2)=:\zeta$ and $z_1,z_2$ are in the same orbit of ${\cal G}_\zeta$. Equivalently, $R$ can be viewed as a subspace of $Z\times_{{\Bbb C}P^1} Z$ over ${\Bbb C}P^1$. A quotient of $Z$ by this relation is a topological space $Z/{\cal G}$ which also has a natural structure of ${\Bbb C}$-ringed space. We define
\begin{definition} An action of ${\cal G}$ on a smooth $Z$ is called {\em regular}, if $Z/{\cal G}$ is a smooth space over ${\Bbb C}P^1$ and the natural projection $\pi:Z\rightarrow Z/{\cal G}$ is a submersion.\label{regular}\end{definition}
A theorem of Godement (\cite{Serre}, Part II, Thm. 3.12.2) gives us a necessary and sufficient condition for regularity:
\begin{proposition} An action of a twistor group ${\cal G}$ on a smooth $Z$ is regular if and only if $R$ is a smooth closed subspace of $Z\times_{{\Bbb C}P^1} Z$ over ${\Bbb C}P^1$. In particular, if the action is free, then it is regular if and only if the quotient $Z/{\cal G}$ is Hausdorff, i.e. $R$ is a smooth closed subspace of $Z\times_{{\Bbb C}P^1} Z$ over ${\Bbb C}P^1$.\label{quotient}\end{proposition}

Suppose now that the action of ${\cal G}$ on $Z$ is free and regular.  Then we shall call $Z$ a {\em principal ${\cal G}$-bundle} over $T:=Z/{\cal G}$. We can classify these as follows:
\begin{proposition} Let $T@>p>> {\Bbb C}P^1$ be a smooth space over ${\Bbb C}P^1$ equipped with a real structure. The set of isomorphism classes of principal ${\cal G}$-bundles over $T$ is in bijective correspondence with elements of $H^1_{\Bbb R}\bigl(T,{\cal O}(p^\ast {\cal G})\bigr)$.\label{principal}\end{proposition}
\begin{pf} Let $P$ be a principal ${\cal G}$-bundle over $T$ and let $\pi_o:p^\ast{\cal G}\rightarrow T$ denote the ``trivial" ${\cal G}$-bundle. Since the projection $\pi:P\rightarrow T$ is a submersion, it admits local sections through every point. Therefore we have ${\cal G}_{p(U_i)}$-equivariant isomorphisms $h_i:\pi^{-1}(U_i)\rightarrow \pi_o^{-1}(U_i)$ for some covering $U_i$ of $T$. For each nonempty intersection $U_i\cap U_j$ we have the transition function $\phi_{ij}=h_jh_i^{-1}$ which is a ${\cal G}$-equivariant (fibrewise) automorphism of $\pi_o^{-1}(U_i\cap U_j)$. Let $G$ denote a particular fiber of $p^\ast{\cal G}$ and let $\phi$ denote $\phi_{ij}$ restricted to this fiber. Thus $\phi:G\rightarrow G$ and $\phi(gx)=g\phi(x)$ for any $x\in G$. Therefore $\phi$ is determined by the value of $\phi(1)$.  One checks that the cocycle conditions for local equivariant automorphisms of $p^\ast{\cal G}$ and for local sections of $p^\ast{\cal G}$ coincide, and this concludes the proof. \end{pf}

\subsection{Orbits and homogeneous spaces} 

The  definition of an orbit of a twistor group as an orbit of a section is quite inadequate. We remark that in the case of an action $\cdot:G\times M\rightarrow M$ of a Lie group $G$ on a manifold $M$, an orbit can be defined in two ways: 1) as the image of a point $m$ under the mapping $\enskip\cdot m:G \rightarrow M$, or 2) as a $G$-homogeneous $G$-invariant submanifold of $M$. The second definition is more suitable in the case of twistor groups.
\begin{definition} Let a twistor group ${\cal G}@>\pi>> {\Bbb C}P^1$ act on a space $Z$ over ${\Bbb C}P^1$ equipped with a real structure $\tau$. An {\em orbit} of ${\cal G}$ is a $\tau$-invariant ${\cal G}$-homogeneous subspace of $Z$.\label{orbit}\end{definition}  

Thus to know the structure of possible orbits we should classify the homogeneous ${\cal G}$-spaces. For example, if ${\cal G}$ is vector bundle (with additive group structure), then ${\cal G}$ acts transitively on any affine bundle $A$ (equipped with an appropriate real structure) such that the linear part of its transition function coincides with the transition function of ${\cal G}$. Thus $A$ is an orbit of ${\cal G}$. 
\par
Firts of all we have
\begin{proposition} Let ${\cal H}$ be a closed twistor subgroup of a twistor group ${\cal G}$. Then ${\cal G}/{\cal H}$ is a smooth homogeneous ${\cal G}$-space.\end{proposition}
\begin{pf} Since the projection ${\cal G}\rightarrow {\Bbb C}P^1$ is a submersion, ${\cal G}$ admits local sections through every point. Let $g(\zeta)$ be such a section over an $U\subset {\Bbb C}P^1$. Consider $\Lie({\cal G})$ and its subalgebra $\Lie({\cal H})$. By taking a smaller $U$, we can assume that $\Lie({\cal G})$ is trivial as a vector bundle over $U$. Therefore there is a subbundle $V$ of $\Lie({\cal G})_{|U}$ complementary to $\Lie({\cal H})$. $V$ is also trivial and we have local section $e_1,\dots,e_m$, which provide a basis of $V_\zeta$ for each $\zeta$.  Then the map $(\zeta,v_1,\dots,v_m)\mapsto g(\zeta)\exp_\zeta\bigl(v_1e_1(\zeta)+\dots+v_me_m(\zeta)\bigr)$ composed with the projection ${\cal G}\rightarrow {\cal G}/{\cal H}$ provides a local coordinate system on ${\cal G}/{\cal H}$. Here $\exp_\zeta$ denotes the fibrewise exponential mapping. The fact that ${\cal G}$ admits a local section through every point implies that so does ${\cal G}/{\cal H}$ and hence the projection ${\cal G}/{\cal H}\rightarrow {\Bbb C}P^1$ is a submersion.\end{pf}

From Proposition \ref{principal} we already know homogeneous ${\cal G}$-spaces on which ${\cal G}$ acts freely: they are given by elements of $H^1_{\Bbb R}({\Bbb C}P^1, {\cal O}({\cal G}))$. In fact, the same argument allows to classify  homogeneous ${\cal G}$-spaces which are locally isomorphic to ${\cal G}/{\cal H}$: 
\begin{proposition} Let ${\cal H}$ be a closed twistor subgroup of ${\cal G}$. The set of isomorphism classes of homogeneous ${\cal G}$-spaces which are locally isomorphic to ${\cal G}/{\cal H}$ is in bijection with elements of $H^1_{\Bbb R}\left({\Bbb C}P^1, {\cal O}\bigl(N({\cal H})/{\cal H}\bigr)\right)$, where $N({\cal H})$ denotes the normalizer of ${\cal H}$.\end{proposition}
\begin{pf} We proceed as in the proof of Proposition \ref{principal} obtaining transition functions $\phi_{ij}$ for such a homogeneous space, which are ${\cal G}$-equivariant (fibrewise) automorphisms of $\pi^{-1}(U_i\cap U_j)$, where $\pi$ denotes the projection ${\cal G}/{\cal H}\rightarrow {\Bbb C}P^1$. Let $G$ and $H$ denote a particular fiber of ${\cal G}$ and ${\cal H}$ and let $\phi$ denote $\phi_{ij}$ restricted to this fiber. Thus $\phi:G/H\rightarrow G/H$ and $\phi(gx)=g\phi(x)$ for any $x\in G/H$. Thus $\phi$ is determined by the value of $\phi(H)$. Suppose that $\phi(H)=pH$. Then, since $\phi(H)=\phi(hH)=h\phi(H)=hpH$, $p\in N(H)$ and such $\phi's$ are in 1-1 correspondence with elements of $N(H)/H$. Once again the cocycle conditions for local equivariant automorphisms of ${\cal G}/{\cal H}$ and for local sections of $N({\cal H})/{\cal H}$ coincide, and this concludes the proof. \end{pf}

In general, there is no reason why a homogeneous ${\cal G}$-space should be locally isomorphic to a fixed ${\cal G}/{\cal H}$. We have, however:
\begin{proposition} Let a twistor group ${\cal G}$ act transitively on a smooth space $W$ in such a way that, for each $\zeta\in {\Bbb C}P^1$, the stabilizer of the action of the fiber $G_\zeta$ on $W_\zeta$ is a normal subgroup of $G_\zeta$. Then there exists a closed normal subgroup ${\cal H}$ of ${\cal G}$ such that $W$ is locally isomorphic to ${\cal G}/{\cal H}$. \end{proposition}
\begin{pf} Since $W$ is smooth, the projection $p$ on ${\Bbb C}P^1$ is a submersion, and so local sections of $W$ exist. Using the transitivity, it follows that locally $p^{-1}(U)$ is ${\cal G}_U$-equivariantly isomorphic to ${\cal G}_U/H(U)$, where $H(U)$ is a subgroup of ${\cal G}_U$. Consider the intersection of two such sets $U_1$ and $U_2$ and proceed as in the proof of the previous proposition. This time, on each fiber, we have a $G$-equivariant diffeomorphism $\phi$ from $G/H_1$ to $G/H_2$. As in the previous proof, $\phi$ is completely determined by its value on $H_1$, say $pH_2$. It follows that $p^{-1}H_1p=H_2$, and since $H_1$ is normal, $H_1=H_2$. Therefore there is a subgroup ${\cal H}$ of ${\cal G}$ whose restriction to each $U$ is $H(U)$.\end{pf}

\section{Negative twistor groups and deformations of hyperk\"ahler structures}

Let a twistor group ${\cal G}$ act regularly (i.e. the quotient is smooth) on a twistor space $Z$ of a hyperk\"ahler manifold $M$ (i.e. $Z/{\cal G}$ is smooth and the natural projection $\pi:Z\rightarrow Z/{\cal G}$ is a submersion). The space $Z/{\cal G}$ over ${\Bbb C}P^1$ has an induced real structure and the pre-image $\pi^{-1}(s)$ of any real section of $Z/{\cal G}$ is a ${\cal G}$-orbit in the sense of definition \ref{orbit}. Therefore we define:
\begin{definition} With the above assumptions, the (real-analytic) space of real sections of $Z/{\cal G}\rightarrow {\Bbb C}P^1$ is called the {\em space of ${\cal G}$-orbits} and is denoted by $M/{\cal G}$. \label{orbitspace}\end{definition}
We have a natural map $\rho:M\rightarrow M/{\cal G}$ obtained by projecting the twistor lines corresponding to points of $M$ to $Z/{\cal G}$. Consider a section $s$ of $Z/{\cal G}$ corresponding to such $\pi(m)$, $m\in M$. Suppose that the action of ${\cal G}$ is locally free. If $L$ denotes $\Lie({\cal G})$, we have an exact sequence of vector bundles:
\begin{equation}0\rightarrow L\rightarrow E\rightarrow T^{V}_s(Z/{\cal G}) \rightarrow 0, \label{sq1}\end{equation}
where ${ V}$ denotes the vertical bundle and $E$ is the normal bundle to the twistor line, i.e. $O(1)\otimes {\Bbb C}^{2n}$. It follows that the normal bundle $N$ to $s$ in $Z/{\cal G}$ splits as the sum of line bundles of degree at least $1$. Therefore $H^1_{\Bbb R}({\Bbb C}P^1,{\cal O}(N))=0$ and the well-known theorem of Kodaira shows that any section of $N$ can be integrated to a deformation the section $s$. This makes a neighbourhood of $\rho(M)$ in $M/{\cal G}$ into a smooth manifold of dimension $\dim M - h^0(L)+h^1(L)$.

We shall now restrict our attention to a special kind of twistor groups. We adopt the following definition
\begin{definition} A twistor Lie algebra ${\cal L}$ is called {\em negative} if it is a direct sum of line bundles of negative degree. A smooth twistor group is called {\em negative} if its Lie algebra is negative. \end{definition}
We also adopt the convention that the identity twistor group $\text{Id}:{\Bbb C}P^1\rightarrow {\Bbb C}P^1$ is also negative.
\par
We have the following simple properties of negative twistor groups which are consequence of Proposition \ref{neg} and the considerations preceding it:
\begin{proposition} \begin{itemize}
\item[(1)] A negative twistor group is nilpotent.
\item[(2)] Any twistor group ${\cal G}$ contains a unique maximal negative subgroup $N({\cal G})$. 
\item[(3)] If $\Lie({\cal G})$ is a direct sum of line bundles of nonpositive degree (e.g. there is a Hamiltonian action of ${\cal G}$ on some twistor space), then $N({\cal G})$ is a normal subgroup of ${\cal G}$ and ${\cal G}/N({\cal G})$ is a trivial group (with some structure group $H$). \hfill $\Box$\end{itemize}\end{proposition}

It is easy to classify smooth twistor groups whose Lie algebra is a fixed negative ${\cal L}$. First of all, there exists a unique twistor group ${\cal G}$ with simply connected fibers and such that $\Lie ({\cal G})={\cal L}$. Indeed, as a manifold ${\cal G}={\cal L}$ and, since every fiber of ${\cal L}$ is a nilpotent Lie algebra, the fibrewise multiplication in ${\cal G}$ is determined by the Campbell-Hausdorff formula. 
\par
By doing things once again fibrewise, we conclude that any other twistor group with the same Lie algebra is of the form ${\cal G}/{\cal D}$ for some (fibrewise) discrete twistor subgroup of the center of ${\cal G}$. 
\par
We shall usually denote negative twistor groups by the letter ${\cal N}$. 
\par
For us, the importance of negative twistor groups follows from the following observation:
\begin{proposition} Let a negative twistor group ${\cal N}$ act regularly a hyperk\"ahler manifold $M$. Then $\rho:M\rightarrow M/{\cal N}$ is an imbedding.\end{proposition}
\begin{pf} To see that $\rho$ is an immersion observe that $d\rho$ at any point of $M$, i.e. at a section of $Z$, is given by the long exact sequence of cohomology induced by \eqref{sq1}. Since $H^0_{\Bbb R}({\Bbb C}P^1,L)=0$, $d\rho$ is injective. Let us show that $\rho$ is injective  $M/{\cal N}$ corresponds to ${\cal N}$-orbits in $M$. The map $\rho$ assigns to a section of the twistor space $Z$, corresponding to a point $m\in M$, its ${\cal N}$-orbit. Thus $\rho$ fails to be an imbedding if an orbit of a section of $Z$ admits more than one section. Since such an orbit $W$ admits a section, it is of the form ${\cal N}/{\cal H}$ for some closed twistor subgroup (not necessarily smooth) ${\cal H}$ of ${\cal N}$. Since $W$ admits two sections, ${\cal N}$ contains two different copies of ${\cal H}$ and so two different sections. Now ${\cal N}$ is of the form  ${\cal G}/{\cal D}$ where ${\cal G}$ is isomorphic to $\Lie({\cal N})$ and repeating the argument implies that ${\cal G}$ and so $\Lie({\cal N})$ contains two sections which is impossible.\end{pf}

Now suppose that the action of ${\cal N}$ on $Z$ is, in addition to being regular, almost free. Then, according to Proposition \ref{principal}, the fibration $Z\rightarrow Z/{\cal N}$ comes from a $\tau$-invariant element of $H^1_{\Bbb R}(Z/{\cal G}, {\cal O}({\cal N}))$. Restricting this cohomology class to sections of $Z/{\cal G}$ gives us a map
\begin{equation} \Lambda:M/{\cal N}\rightarrow H^1_{\Bbb R}({\Bbb C}P^1, {\cal O}({\cal N})).\label{Lambda}\end{equation}
At this point a remark about the structure of $H^1_{\Bbb R}({\Bbb C}P^1, {\cal O}({\cal N}))$ is in order. It is not a group, unless ${\cal N}$ is abelian. It does have, however, a prefered element $1$ (corresponding to ${\cal N}$). In addition, it has a natural structure of a smooth manifold, with charts diffeomorphic to $H^1_{\Bbb R}\bigl({\Bbb C}P^1, \Lie({\cal N})\bigr)$. We observe that the map $\Lambda$ is a smooth, with the differential defined as follows. Let $A$ be an ${\cal N}$-orbit in $Z$, corresponding to an element $\pi(A)$ of $M/{\cal N}$. Then we have an exact sequence of vector bundles
\begin{equation} 0\rightarrow (T^VA)/{\cal N} \rightarrow \bigl(T^V_A Z\bigr)/{\cal N} \rightarrow T^V_{\pi(A)}\bigl(Z/{\cal N}\bigr) \rightarrow 0,\label{sq2}\end{equation}
where the action of ${\cal N}$ on $T^VZ$ is the tangent action along the fibers. The differential of $\Lambda$ at $\pi(A)$ is then the induced map \begin{equation} d\Lambda_{\pi(A)}:H^0_{\Bbb R}\left(\pi(A), T^V_{\pi(A)}\bigl(Z/{\cal N}\bigr)\right) \rightarrow  H^1_{\Bbb R}\bigl(A/{\cal N}, (T^VA)/{\cal N}\bigr)\simeq H^1_{\Bbb R}\bigl({\Bbb C}P^1,\Lie({\cal N})\bigr).\label{sq3}\end{equation}
 
 We observe that $\Lambda^{-1}(1)$ corresponds to ${\cal N}$-orbits possesing a section and so, from the previous proposition, to $M$. In general, for any $\lambda$, $\Lambda^{-1}(\lambda)$ parameterizes orbits of the fixed type $\lambda$. We claim that  $M_\lambda:=\Lambda^{-1}(\lambda)$ carries a natural hyperk\"ahler structure, which should be viewed as a deformation of the hyperk\"ahler structure of $M$. More precisely:
\begin{theorem} Let a negative twistor group ${\cal N}$ act regularly, almost freely, and symplectically on a hyperk\"ahler manifold $M^{4n}$. Then there exists a smooth neighbourhood $U$ of $M$ in $M/{\cal N}$ such that $\Lambda$ is a submersion on $U$ and, for any $\lambda\in  H^1_{\Bbb R}({\Bbb C}P^1, {\cal O}({\cal N}))$, $M_\lambda:=\Lambda^{-1}(\lambda)\cap U$ carries a natural hyperk\"ahler structure. Furthermore, with respect to each complex structure, $M_\lambda$ is isomorphic, as a complex symplectic manifold, to an open subset of $M$.
\label{deformation} \end{theorem}
\begin{remark} Completely analogous result holds for hypercomplex manifolds. \end{remark}
\begin{pf} We consider the vector bundle $F=\bigl(T^V Z\bigr)/{\cal N}$ on $Z/{\cal N}$. Over a section obtained by projecting a twistor line $s$ in $Z$, this bundle is just the normal bundle of $s$, and so $O(1)\otimes{\Bbb C}^{2n}$. By standard semi-continuity theorems, $F$ is  $O(1)\otimes{\Bbb C}^{2n}$ when restricted to neighbouring sections, i.e. to a neighbourhood $U$  of $\rho(M)\simeq M$ in $M/{\cal N}$. Then \eqref{sq2} and \eqref{sq3} show that $\Lambda$ is a submersion on $U$. Thus $M_\lambda=\Lambda^{-1}(\lambda)\cap U$ is a submanifold of $M/{\cal N}$. The tangent space $T_pM_{\lambda}$ is the space of real sections of $F_{|s(p)}$, where $s(p)$ is the section of $Z/{\cal N}$ given by $p$. This is the same as $\bigl(T^V_A Z\bigr)/{\cal N}$ where $A$ is the ${\cal N}$-orbit in $Z$, whose projection is $s(p)$. We have an  $O(2)$-valued complex-symplectic form on the fibers of $\bigl(T^V_A Z\bigr)/{\cal N}$, given by $\tilde{\omega}([a],[b])=\omega(a,b)$, where  $\omega$ is the given form on $Z$ and the representatives $a,b$ are tangent to the same point of $A$. Since ${\cal N}$ acts symplectically, this does not depend on the choice of point in $A$. We notice that on each fiber over ${\Bbb C}P^1$ this is canonically isomorphic to $\omega$ on this fiber. In particular $\tilde{\omega}$ is nondegenerate and closed. Now, as in the proof of Theorem   in \cite{HKLR}, we obtain a hyperhermitian structure on $M_\lambda$. The above isomorphisms on each fiber give us local isomorphisms of complex structures (essentially $(Z_{\zeta}\times N_\zeta)/N_\zeta\simeq Z_\zeta$) proving their integrability and proving the theorem.\end{pf}

The above proof allows us to identify the twistor space of $M_\lambda$. Let $W_\lambda$ be a principal ${\cal N}$-bundle over ${\Bbb C}P^1$ corresponding to a $\lambda\in H^1_{\Bbb R}({\Bbb C}P^1, {\cal O}({\cal N})) $. Let $Z$ be the twistor space of $M$ and consider the diagonal action of ${\cal N}$ on $Z\rightarrow W_\lambda$. Then $Z_{\lambda}=\bigl(Z\times_{{\Bbb C}P^1} W_\lambda\bigr)/{\cal N}$ is the twistor space of $M_\lambda$ (i.e. $M_\lambda$ is the family of sections of $Z_\lambda$ with correct normal bundle). We observe that ${\cal N}$ does not necessarily act on $Z_{\lambda}$. It acts only if ${\cal N}$ is abelian. In general case, we obtain an action of another negative twistor group ${\cal N}_\lambda$, locally isomorphic to ${\cal N}$ and obtained by gluing pieces of ${\cal N}$ by inner automorphisms corresponding to local sections of ${\cal N}$ determined by $\lambda\in H^1_{\Bbb R}({\Bbb C}P^1, {\cal O}({\cal N}))$ (for $\lambda$ close to $1$, the Lie algebra of ${\cal N}_\lambda$ must be negative).

\section{Twistor quotients\label{Quotient}}

We now wish to associate a ``quotient" to a hyperk\"ahler manifold with a twistor group action. Essentially, this quotient is formed by taking the complex symplectic quotients along the fibers of the twistor space. 
\par
Let therefore a twistor group ${\cal G}@>\pi>> {\Bbb C}P^1$ act  on the twistor space $Z@>p>> {\Bbb C}P^1$ of a hyperk\"ahler manifold $M$. We suppose that this action is Hamiltonian with the moment map $\mu:Z\rightarrow {\cal L}^\ast\otimes {\cal O}(2)$. Here ${\cal L}=\Lie{\cal G}$.
\par
Let $s$ be a twistor line in $Z$. Then $\mu\circ s$ is a real section of ${\cal L}^\ast\otimes {\cal O}(2)$. Let $S=(\mu\circ s)({\Bbb C}P^1)$ and suppose that
the fibrewise quotient of $\mu^{-1}(S)$ by $\text{Stab}(\mu\circ s)$ (stabilizer of coadjoint action) is a manifold (fibering over ${\Bbb C}P^1$), which we denote by $Z_{\text{red}}$. It inherits the real structure, the twisted complex-symplectic form along the fibers and a real section $\bar{s}$, induced by $s$. Thus, if $Z_{\text{red}}$ contains a real section (e.g. $\bar{s}$) whose normal budle is the sum of line bundles of degree $1$, then $Z_{\text{red}}$ is a twistor space of a pseudo-hyperk\"ahler manifold, which we denote by $M//{\cal G}$. We shall call this construction the  {\em ``twistor quotient}. If $M$ has dimension $4n$ and the complex dimension of the fiber of $\Lie{\cal G}$ is $m$, then $M//{\cal G}$ has dimension $4n-4m$.
\par
What we need then are conditions which guarantee that  $Z_{\text{red}}$ has sections with correct normal budle. First of all, if the action of ${\cal G}$ is locally free at $s$, then the normal bundle of $\bar{s}$ in $Z_{\text{red}}$ is $L^\perp/(L\cap L^\perp)$, where $L$ is the subbundle of the normal bundle of $s$ in $Z$ generated by the Lie algebra of ${\cal G}$ and the orthogonal is taken with respect to the twisted symplectic form. Let us make the following definition.
\begin{definition} Let $s$ be a twistor section of a twistor space $Z@>p>> {\Bbb C}P^1$ of a hyperk\"ahler manifold $M$ on which there is a locally free Hamiltonian action of a twistor group ${\cal G}$.  We shall say that $s$ is {\em ${\cal G}$-admissible} if $L^\perp/(L\cap L^\perp)$ is the sum of line bundles of degree $1$.\end{definition}

Thus, if $s$ is ${\cal G}$-admissible and $Z_{\text{red}}$ is a manifold in a neighbourhood of $s$, then $Z_{\text{red}}$ is a twistor space of a pseudo-hyperk\"ahler manifold. We now have:

\begin{proposition} Let ${\cal H}$ be a twistor subgroup of a twistor group ${\cal G}$ such that $\Lie({\cal G})/\Lie({\cal H})$ is the sum of line bundles of degree $1$. Suppose that we have a locally free Hamiltonian action of ${\cal G}$ on a $Z@>p>> {\Bbb C}P^1$ with a moment map $\mu$. Then any ${\cal G}$-admissible twistor section $s$ of $Z$ such that $\mu\circ s$ is ${\cal G}$-invariant, is ${\cal H}$-admissible.\label{criterion}\end{proposition}
\begin{pf} As above, let $L$ denote the subbundle of the normal bundle of $s$
generated by the action of ${\cal G}$. Since ${\cal G}$ acts locally freely at $s$, $L\simeq \Lie({\cal G})$ as vector bundles. Furthermore, since $\mu\circ s$ is ${\cal G}$-invariant, $L\subset L^\perp$. Let $P$ denote the subbundle of $L$ generated by ${\cal H}$. We have to show that $P^\perp/P$ is the sum of line bundles of degree $1$. We observe that it is enough to show that $H^1\bigl((P^\perp/P)^\ast\bigr)=0$. Indeed, this implies that $P^\perp/P$ is the sum of line bundles of degree at most $1$, and since we also have the isomorphism  $P^\perp/P \simeq (P^\perp/P)^\ast\otimes {\cal O}(2)$ given by the $\omega$, all line bundle summands in  $P^\perp/P$ have degree $1$. 
\par
To show that $H^1\bigl((P^\perp/p)^\ast\bigr)=0$ it is sufficient to show that the map $H^0\bigl((P^\perp)^\ast\bigr)\rightarrow H^0(P^\ast)$ is surjective (as $P^\perp$ is a subbundle of the normal bundle of $s$ - sum of line bundles of degree $1$ - therefore $H^1\bigl((P^\perp)^\ast\bigr)=0$).  We have the following embeddings of vector bundles
$$ P\hookrightarrow L \hookrightarrow L^\perp \hookrightarrow P^\perp. $$
We shall show that the dual of each of these maps is surjective on $H^0$ by showing that $H^1$ of each quotient vanishes. For the first one, $H^1\bigl(L/P)^\ast \bigr)=0$ by our assumption. For the middle map this follows from the fact that $L^\perp/L$ is the sum of ${\cal O}(1)$'s (by assumption, $s$ is ${\cal G}$-admissible). For the last one, we have to show that $H^1\bigl((P^\perp/L^\perp)^\ast\bigr)=0$. The form $\omega$ and Serre duality show that this cohomology group is the same as $H^0\bigl((P/L)^\ast \bigr)$ which again vanishes by our assumption. \end{pf}
There are two cases, when a twistor section $s$ is automatically ${\cal G}$-admissible: 1) if ${\cal G}$ is trivial twistor group $\underline{G}$; and 2) if $L\simeq \Lie({\cal G})$ is a Lagrangian subbundle of the normal bundle of $s$. This second condition holds, e.g., in the case of the generalized Legendre transform. We make several other remarks:   

\begin{remark} A necessary condition for $\Lie(\underline{G})/\Lie({\cal H})$ to be the sum of ${\cal O}(1)$'s is that, as a vector bundle, $\Lie(\underline{H})=\sum {\cal O}(-p_i)$ with $\sum p_i=d$, where $d$ is the fiber codimension of ${\cal H}$ in $\underline{G}$.\end{remark}

\begin{remark} The sufficient condition of this proposition is particularly useful when dealing with abelian twistor groups. If both $\underline{G}$ and ${\cal H}$ are abelian (e.g. vector bundles), and the numerical condition of the previous remark is satisfied, then a generic embedding of ${\cal H}$ into  $\underline{G}$ gives a twistor quotient. Thus, for example, if there is a locally free $3$-Hamiltonian action of ${\Bbb R}^3$ (effective or not) on a hyperk\"ahler manifold $M$ which extends to an action of ${\Bbb C}^3$ with respect to each complex structure, then a generic embedding of ${\cal O}(-2)$ (compatible with real structures) into $\underline{{\Bbb C}^3}$ satisfies the condition of Proposition \ref{criterion}. \end{remark}

There also is a simple necessary condition in the setting of  Proposition \ref{criterion}. Namely, since (in the notation of the proof of that proposition) $L/P\hookrightarrow P^\perp/P$, we need that $\Lie({\cal G})/\Lie({\cal H})$ is the sum of line bundles of degree {\em at most} $1$.  Thus we shall not find twistor quotients by  ${\cal O}(-4)$ embedded into $\underline{{\Bbb C}^3}$.

Let us turn to examples. 
\begin{example} Santa Cruz \cite{SC} constructed twistor spaces of hyperk\"ahler metrics on coadjoint orbits of complex semisimple Lie groups (see also \cite{AG}). He associates such a metric to any real section (spectral curve) $s$ of ${\frak  g}\otimes O(2)$, whose fibrewise stabilizers have constant dimension. Here ${\frak  g}$ is the Lie algebra of a Lie group $G$. His construction can be interpreted as a twistor quotient of the hyperk\"ahler metric on $T^\ast G$ (cf. \cite{BielJLMS}) by the trivial twistor group $\underline{G}$, where the level set of the moment map is chosen to be $s$. In other words the resulting twistor space is $\mu^{-1}(s)/\text{Stab}(s)$, i.e. the fibrewise complex-symplectic quotient of the twistor space of $T^\ast G$ by $G$. \end{example}
\begin{example} Many interesting metrics can be constructed as twistor quotients by the group ${\cal N}$ defined in Example \ref{borel}. Thus whenever we have  an effective triholomorphic and isometric action of a compact Lie group $G$ on a hyperk\"ahler manifold $M$ we can form a twistor quotient of $M$ by ${\cal N}$. This is a reinterpretation of the construction given in \cite{BielJLMS}. 
\par 
In particular, the natural hyperk\"ahler metric on the moduli space of $SU(2)$-monopoles of charge $k$ can be obtained as such a quotient of $T^\ast Gl(k, {\Bbb C})$. Also the ALE spaces can be obtained as such quotients of coadjoint orbits with Kronheimer's metric \cite{BielAGAG}.\end{example}
\par
It is possible to know the metric on the twistor quotient of $M$, if we know the metric on the deformations $M_\lambda$ of Theorem \ref{deformation}. First of all, since a twistor group ${\cal G}$, by which we quotient, admits a chain of subgroups ${\cal G}={\cal G}_1\subset \ldots \subset {\cal G}_k$, such that each subgroup is normal in the previous one and ${\cal G}_i/ {\cal G}_{i+1}$ is abelian for $i<k-1$ and trivial for $i=k-1$, a twistor quotient by an arbitrary group reduces to twistor quotients by abelian twistor groups and to hyperk\"ahler quotients. We shall, therefore, assume for the remainder of the section that ${\cal G}$ is abelian. In this case the moment map $\mu:Z\rightarrow \Lie({\cal G})\otimes O(2)$ descends to $Z/{\cal G}$.
\par
Let us choose local complex coordinates  $u_1,\dots, u_n, z_1,\dots, z_n$ in a fiber $Z_\zeta$ of $Z$, so that $u_1,\ldots, u_k$ correspond to ${\cal G}_\zeta$ and $ z_1,\ldots,z_k$ give us the complex moment map for ${\cal G}_\zeta$. The remaining coordinates give complex coordinates on $M//{\cal G}$ 
\par
Since $\mu$ descends to $Z/{\cal G}$, we have a corresponding map
$$\Phi:M//{\cal G}\rightarrow \Gamma(\Lie({\cal G})\otimes O(2).$$
$M//{\cal G}$ can be identified with $\Phi^{-1}(v_0)$, for some $v_0\in V=\Gamma(\Lie({\cal G})\otimes O(2)$. The coordinates on $M/{\cal G}$ are $u_{k+1},\dots,u_n,z_1,\dots,z_n$ and the remaining (apart from $z_1,\dots,z_k$) coordinates of $V$. Now, on each deformation $\Lambda^{-1}(\lambda)\subset M/{\cal G}$ of Theorem \ref{deformation} we have a K\"ahler potential 
$$K_\lambda(u_1,\dots, u_n, z_1,\dots, z_n)=K(u_{k+1},\dots, u_n, z_{k+1},\dots, z_n, v)$$
where $v$ varies over $ V$. The K\"ahler potential for  $M//{\cal G}$ is then
$$\bar{K}(u_{k+1},\dots, u_n, z_{k+1},\dots, z_n)=(u_{k+1},\dots, u_n, z_{k+1},\dots, z_n, v_0).$$

\section{The generalized Legendre transform}
Linndstr\"om and Ro\v{c}ek \cite{LR} found several constructions of hyperk\"ahler metrics, in particular two based on the Legendre transform. The simpler one produces precisely hyperk\"ahler metrics in $4n$ dimensions with a local tri-Hamiltonian (hence isometric) action of ${\Bbb R}^n$. The second one, the generalized Legendre transform (GLT), produces metrics which generically don't have triholomorphic vector fields. 
\par In the simplest case of $4$-dimensional metrics, such a metric is associated to a real-valued function $F$ on ${\Bbb R}^{2k+1}$, $k\geq 2$,
 with coordinates $w_0,\ldots, w_{2k}\in {\Bbb C}$, $w_{2k-i}=(-1)^{k+i}$ which satisfies the system of linear PDE's:
$$ F_{w_iw_j}=F_{w_{i+s}w_{j-s}}$$
for all $i,j,s$. The hyperk\"ahler metric lives on the submanifold of ${\Bbb R}^{2k+1}$, defined by the equations $F_{w_i}=0$, for $2\leq i\leq 2k-2$. 
An example of a metric which can be constructed using the GLT is the Atiyah-Hitchin metric \cite{I-R} or other $SU(2)$-monopole metrics \cite{Hou}.
\par
In \cite{I-R}, Ivanov and Ro\v{c}ek gave an interpretation of metrics constructed via GLT in terms of twistor spaces. They show that the twistor space $Z$ of such a manifold $M$ has a projection $p$ onto $O(2k)$ (or at least an open subset of it, invariant under the real structure).  Moreover the kernel of $dp$ is a Lagrangian subbundle of $T^VZ$. We can interpret this as saying that there is a local action of the twistor group $O(-2k+2)$ on $Z$. The projection onto $O(2k)$ is the moment map, and $O(2k)$ can be identified (if the fibers of $p$ are connected) with $Z/O(-2k+2)$. The vector space ${\Bbb R}^{2k+1}$ is $M/O(-2k+2)$ and the equations  $F_{w_i}=0$, for $2\leq i\leq 2k-2$, which determine $M$, are equivalent to setting the $\Lambda$ of \eqref{Lambda} equal to zero. 
\par
A similar interpretation holds for higher dimensional hyperk\"ahler metrics constructed via the generalized Legendre transform. This construction produces $4n$-dimensional metrics which admit a local Hamiltonian action of an $n$-dimensional abelian twistor group.

\begin{ack} This work has been supported by EPSRC's Advanced Research Fellowship, which is gratefully acknowledged. I also thank Martin Ro\v{c}ek for useful discussions. \end{ack}

\addtolength{\textheight}{1cm}

\end{document}